\documentclass[12pt,draftcls,onecolumn]{IEEEtran}

\newtheorem{CC}{Corollary}
\newtheorem{PP}{Proposition}

\newcommand\norm[1]{\left\lVert#1\right\rVert}

\ifCLASSINFOpdf
\else
\fi

\ifCLASSOPTIONcompsoc
\usepackage[nocompress]{cite}
\else
\usepackage{cite}
\fi

\usepackage{amsmath}
\usepackage{amssymb}
\usepackage{graphicx}
\usepackage{epstopdf}
\usepackage{mathtools}
\usepackage{xcolor}

\begin{document}

\title{A Matrix-Valued Inner Product for Matrix-Valued Signals and Matrix-Valued Lattices}


\author{
        Xiang-Gen Xia 
\vspace{-0.5in}

\thanks{
X.-G. Xia is with the Department of Electrical and Computer Engineering, University of Delaware, Newark, DE 19716, USA (e-mail: xxia@ece.udel.edu).}

}

\date{}

\maketitle


\begin{abstract}
A matrix-valued inner product was proposed before to construct 
orthonormal matrix-valued wavelets for matrix-valued signals.
It introduces a weaker orthogonality for matrix-valued signals
than the orthogonality of all components in a matrix that is 
commonly used in orthogonal multiwavelet constructions. 
With the weaker orthogonality, it is easier to construct orthonormal 
matrix-valued wavelets. 
In this paper, we re-study the matrix-valued inner product more from 
the inner product viewpoint that is more fundamental and 
propose a new but equivalent norm for matrix-valued signals.
 We show that although it is not scalar-valued, 
it maintains most of the scalar-valued inner product properties.
We introduce a new linear independence concept for matrix-valued
signals and present  some related properties. We then present 
the Gram-Schmidt orthonormalization procedure
for a set of linearly independent matrix-valued signals. 
Finally we  define matrix-valued lattices, 
where the newly introduced Gram-Schmidt orthogonalization might be applied. 
\\
\\

\end{abstract}



\section{Introduction}\label{sec1}
Matrix-valued (vector-valued) signals 
are everywhere these days, such as, videos, 
multi-spectral images, signals from multiarray multisensors,
and high dimensional data. 
For these signals, there are correlations not only over the time
but also across their matrix components. 
How to efficiently represent and process them
plays an important and fundamental role in data science. 

In \cite{xia1,xia2}, a matrix-valued inner product was introduced 
for matrix-valued (or vector-valued) signals.
It leads to a weaker orthogonality (called Orthogonality B) 
than the component-wise orthogonality (called Orthogonality A) 
in orthogonal multiwavelet constructions \cite{strang, plonka}.
The weaker orthogonality, i.e., Orthogonality B, 
provides an easier sufficient condition to construct  orthonormal multiwavelets with Orthogonality B 
 than the necessary and sufficient 
condition \cite{plonka} to construct  orthonormal multiwavelets
with Orthogonality A. 
A connection between multiwavelets and matrix-valued/vector-valued
wavelets can be found in \cite{xia2}. After the works in \cite{xia1, xia2}, there have been
many studies on matrix-valued/vector-valued wavelets 
for matrix-valued/vector-valued signals 
in the literature, see, for example, \cite{wave1}-\cite{wave6}. 

On the other hand, the Orthogonality B induced 
from the matrix-valued inner product is stronger than 
the orthogonality, which is called  Orthogonality C here, 
  induced from the commonly used scalar-valued
inner product for matrices. 
It is because, as we shall see later, Orthogonality B is basically 
the orthogonality between all row vectors of two  matrix-valued signals,
while  Orthogonality C is the orthogonality of two long vectors of 
concatenated row vectors of two matrix-valued signals. 
It was proved in \cite{xia2} that 
 Orthogonality B or the matrix-valued inner product 
 is  able to  completely decorrelate 
matrix-valued signals not only in time domain but also 
across the components inside matrix in the sense that it provides 
a complete   Karhunen-Lo\`{e}ve  expansion for matrix-valued 
signals,  while Orthogonality A or Orthogonality C  may not do so.
In other words, Orthogonality B induced from the 
matrix-valued inner product  is the proper orthogonality for 
matrix-valued signals.
This also means that the matrix-valued 
inner product is needed to study the decorrelation of matrix-valued
signals and a conventional scalar-valued inner product may not be enough.

Since the main goal in \cite{xia1, xia2} was to construct orthonormal 
matrix-valued (vector-valued) wavelets, not much about the inner 
product or the orthogonality itself, which is more fundamental,
was studied.  
In this paper, we study more properties 
on the matrix-valued inner product and its induced 
Orthogonality B  for matrix-valued 
signal space proposed in \cite{xia1,xia2}.
We first define a different norm for matrix-valued signals than that defined
in \cite{xia2} and prove that these two norms are equivalent. 
The norm defined in this paper is consistent with the matrix-valued 
inner product similar to that for a  scalar-valued inner product. 
We introduce a new linear independence concept for matrix-valued
signals and present some related properties. We then present 
the Gram-Schmidt orthonormalization procedure
for a set of linearly independent matrix-valued signals.
We finally define matrix-valued lattices, 
where the newly introduced Gram-Schmidt orthogonalization might be applied. 
Due to the noncommuntativity of matrix multiplications, these concepts and 
properties for matrix-valued signals and/or inner product 
are not straightforward extensions of the conventional ones
for scalar-valued signals and/or inner product. 

The remainder of this paper is organized as follows. 
In Section \ref{sec2}, we introduce matrix-valued signal space, 
 matrix-valued inner product, and define a new norm for matrix-valued
signals. We present some simple properties for the matrix-valued 
inner product and prove that the new norm proposed in this paper
is equivalent to that used in \cite{xia2}. 
In Section \ref{sec3}, we first introduce the concepts of 
degenerate and nondegenerate 
matrix-valued signals and then introduce the concept 
of  linear independence 
for matrix-valued signals. The newly introduced linear 
independence is different from but consistent with the conventional 
one for vectors. 
We also present some interesting properties on the
linear independence and the orthogonality. 
We finally present the Gram-Schmidt orthonormalization
procedure for a set of linearly independent 
matrix-valued signals, which has the similar form 
as the conventional one for vectors but not a straightforward
generalization due to the noncommutativity of matrix multiplications
and the matrix-valued inner product used in the procedure. 
In Section \ref{sec7}, we define matrix-valued lattices.
In Section \ref{sec8}, we conclude this paper. 

\section{Matrix-Valued Signal Space and Matrix-Valued Inner Product}\label{sec2}

We first introduce matrix-valued signal space studied in \cite{xia1,xia2}.
Let ${\mathbb{C}}^{N\times N}$ denote
all  $N\times N$ matrices of complex-valued entries, 
and for $-\infty \leq a<b\leq \infty$, let $L^2(a,b)$
denote all the finite energy signals in the interval $(a, b)$ and 
\begin{equation}\label{2}
L^2(a,b;{\mathbb{C}}^{N\times N})\stackrel{\Delta}{=}
\left\{{\bf f}(t)=
\left( \begin{array}{cccc}
f_{11}(t) & f_{12}(t) & \cdots & f_{1N}(t)\\
f_{21}(t) & f_{22}(t) & \cdots & f_{2N}(t)\\
          &           & \cdots & \\
f_{N1}(t) & f_{N2}(t) & \cdots & f_{NN}(t)
\end{array}\right):\;\; f_{kl}\in L^2(a,b), 1\leq k,l\leq N\right\}.
\end{equation}
We call $L^2(a,b;{\mathbb{C}}^{N\times N})$ a matrix-valued signal space 
and $\mathbf{f}(t)\in L^2(a,b;{\mathbb{C}}^{N\times N})$, or
simply  $\mathbf{f}\in L^2(a,b;{\mathbb{C}}^{N\times N})$,  a matrix-valued signal.

For any $A\in \mathbb{C}^{N\times N}$ and
${\bf f}\in L^2(a,b;\mathbb{C}^{N\times N})$, the products
$A{\bf f},\,{\bf f}A \in L^2(a,b; \mathbb{C}^{N\times N})$. 
This implies that the  matrix-valued signal space
$L^2(a,b;\mathbb{C}^{N\times N})$ is defined over $\mathbb{C}^{N\times N}$
and not  simply over $\mathbb{C}$.
For $\mathbf{f}\in L^2(a,b; {\mathbb{C}}^{N\times N})$, its integration $\int_a^b \mathbf{f}(t) dt$ is defined by the integrations of its components,
i.e., $\int_a^b \mathbf{f}(t) dt= \left(\int_a^b f_{kl}(t)dt\right)$.

Let $\| \cdot \|_M$ denote a matrix norm on $\mathbb{C}^{N\times N}$,
for example, the normalized Frobenius norm for a matrix, 
$$
\|A\|_M =\|A\|_F= \left( \frac{1}{N} \sum_{k,l=1}^{N}
|A_{kl}|^2 \right)^{1/2},
$$
where $A=(A_{kl})$.
For each ${\bf f}\in L^2(a,b;\mathbb{C}^{N\times N})$, let $\| {\bf f}\|_M$
denote the norm of ${\bf f}$ associated with the matrix norm
$\| \cdot \|_M$ as
\begin{equation}\label{3}
\|{\bf f}\|_M \stackrel{\Delta}{=} 
\norm{ \int_a^b{\bf f}(t){\bf f}^{\dagger}(t)dt}_M^{1/2},
\end{equation}
where $^{\dagger}$ denotes the complex conjugate transpose.
Note that the norm $\|{\bf f}\|$ 
 of ${\bf f}$ defined in \cite{xia2} has the following 
form 
\begin{equation}\label{4}
\|{\bf f}\|\stackrel{\Delta}{=}\left( \int_a^b\|{\bf f}(t)\|_M^2dt \right)^{1/2},
\end{equation}
where $\|{\bf f}(t)\|_M$ is the matrix norm of matrix  ${\bf f}(t)$ for a fixed $t$. 
We will show later that the above two norms $\|{\bf f}\|$ and
$\|{\bf f}\|_M$  are equivalent
in the sense that there exist two positive constants $C_1>0$ and $C_2>0$
such that
\begin{equation}\label{5}
C_1\|{\bf f}\| \leq \|{\bf f}\|_M \leq C_2 \|{\bf f}\|, 
\mbox{ for any } {\bf f}\in L^2(a,b; {\mathbb{C}}^{N\times N}).
\end{equation}

We next define matrix-valued inner product for matrix-valued signals
in $L^2(a,b; {\mathbb{C}}^{N\times N})$. 
For two matrix-valued signals
${\bf f}, {\bf g}\in L^2(a,b; {\mathbb{C}}^{N\times N})$,
their {\em matrix-valued inner product} (or simply inner product) 
$\langle {\bf f},{\bf g}\rangle$ is defined as the integration
of the matrix product ${\bf f}(t){\bf g}^{\dagger}(t)$, i.e.,
\begin{equation}\label{6}
\langle {\bf f},{\bf g}\rangle\stackrel{\Delta}{=}\int_a^b {\bf f}(t){\bf g}^{\dagger}(t)dt.
\end{equation}

With the definition (\ref{6}),  most properties 
of the conventional scalar-valued inner product hold for the
above matrix-valued inner product. For instance, 
the following properties of the matrix-valued inner product  are clear: 
\begin{itemize}
\item[(i)]  $\langle {\bf f},{\bf g} \rangle=\langle {\bf g}, {\bf f}\rangle^{\dagger}$.
\item[(ii)]  $\langle {\bf f},{\bf f} \rangle=0$ if and only if ${\bf f}=0$.

\item[(iii)] $\|{\bf f}\|_M =  \| \langle {\bf f},{\bf f} \rangle \|_M^{1/2}$. 

\item[(iv)] For any $A, B\in {\mathbb{C}}^{N\times N}$, 
$\langle A {\bf f}, B {\bf g}\rangle = A \langle {\bf f}, 
{\bf g}\rangle B^{\dagger}$. 
\end{itemize}
Note that Property (iii) may not hold for 
 the norm (\ref{4}) used in \cite{xia2}. 

Two matrix-valued signals ${\bf f}$ and ${\bf g}$
in $L^2(a,b; {\mathbb{C}}^{N\times N})$ are called {\em orthogonal} 
if $\langle {\bf f},{\bf g} \rangle=0$. 
A set of matrix-valued signals is called an {\em orthogonal set} if
any two distinct matrix-valued signals in the set are orthogonal. 
A sequence $\mathbf{\Phi}_k(t)\in  L^2(a,b; {\mathbb{C}}^{N\times N})$,
 $k\in {\mathbb{Z}}$,
is called an {\em orthonormal set} in $L^2(a,b; {\mathbb{C}}^{N\times N})$ if
\begin{equation}\label{7}
\langle \mathbf{\Phi}_k,\mathbf{\Phi}_l\rangle=\delta(k-l)I_N,\,\,\,k,l\in {\mathbb{Z}},
\end{equation}
where $\delta(k)=1$ when $k=0$ and $\delta(k)=0$ when $k\neq0$, $I_N$ is the $N\times N$ identity matrix.
Due to (i) above, the orthogonality/orthonormality 
between ${\bf f}$ and ${\bf g}$ is commutative, i.e., 
if ${\bf f}$ and ${\bf g}$ are orthogonal/orthonormal, then
${\bf g}$ and ${\bf f}$ are orthogonal/orthonormal too.

 A sequence $\mathbf{\Phi}_k(t)\in  L^2(a,b;{\mathbb{C}}^{N\times N})$, $k\in {\mathbb{Z}}$,
is called an {\em orthonormal basis} for 
$L^2(a,b; {\mathbb{C}}^{N\times N})$ if it satisfies (\ref{7}), and moreover, for any $\mathbf{f}\in L^2(a,b; {\mathbb{C}}^{N\times N})$ there exists a sequence of $N\times N$ constant matrices 
$F_k$, $k\in\mathbb{Z}$,  such that
\begin{equation}\label{8}
\mathbf{f}(t)=\sum_{k\in\mathbb{Z}}F_k\mathbf{\Phi}_k(t),\,\, \mbox{ for }
t\in [a,b],
\end{equation}
or simply
$$
\mathbf{f}=\sum_{k\in\mathbb{Z}}F_k\mathbf{\Phi}_k,
$$
where $F_k=\langle \mathbf{f},\mathbf{\Phi}_k\rangle$, the multiplication $F_k\mathbf{\Phi}_k(t)$
for each fixed $t$ is the $N\times N$ matrix multiplication, and the convergence for the infinite summation is in the sense
of the norm  $\| \cdot \|_M$ defined by (\ref{3}) for the matrix-valued signal space. The corresponding Parseval equality is
\begin{equation}\label{9}
\langle \mathbf{f},\mathbf{f}\rangle=\sum_{k\in\mathbb{Z}}F_kF^{\dagger}_k.
\end{equation}
With the norm $\| \cdot \|_M$ in (\ref{3}), it is clear that for any element 
$\mathbf{\Phi}_k$ in an orthonormal set 
in $L^2(a,b; {\mathbb{C}}^{N\times N})$, 
we have $\|\mathbf{\Phi}_k\|_M=1$, which is consistent 
with the conventional relationship between vector norm and 
vector inner product. 
However, this property may not hold for the norm $\| \cdot \|$ 
in (\ref{4}) used in \cite{xia2}. We refer to \cite{xia2}
for the  Karhunen-Lo\`{e}ve  expansion with an orthonormal 
basis for random processes of matrix-valued signals. 

We next show the equivalence (\ref{5}) 
of the two norms $\| \cdot \|_M$ in (\ref{3})
and $\| \cdot \|$ in (\ref{4}). 

\begin{PP}
\label{p0}
The norms $\|\cdot \|_M$ in (\ref{3}) and 
$\| \cdot \|$ in (\ref{4}) are equivalent.
\end{PP}

{\bf Proof}: 
It is known that all matrix norms 
for constant matrices are equivalent. Hence, we show (\ref{5}) only for the Frobenius norm,
i.e., $\| \cdot \|_M=\| \cdot \|_F$.  In this case, the following is
an elementary proof. 
\begin{eqnarray*}
\|{\bf f}\|_M^4 & = & \frac{1}{N} \sum_{k,l=1}^N
\left| \int_{a}^b \sum_{m=1}^N f_{km}(t) f_{lm}^*(t)dt \right|^2\\
      & = & \frac{1}{N} \sum_{k=1}^N \left| \sum_{m=1}^N \int_a^b |f_{km}(t)|^2dt \right|^2
            + \frac{1}{N} \sum_{k\neq l =1}^N \left| \sum_{m=1}^N \int_a^b
           f_{km}(t) f_{lm}^*(t)dt \right|^2\\
      & \geq & \frac{1}{N}  \sum_{k=1}^N \left| \sum_{m=1}^N \int_a^b |f_{km}(t)|^2dt \right|^2\\
      & \geq & \frac{1}{N^2} \left( \int_a^b \sum_{k,l=1}^N |f_{kl}(t)|^2dt\right)^2\\
      & = & \frac{1}{N^2}  \| {\bf f}\|^4.
\end{eqnarray*}
Thus, we have 
\begin{equation}\label{10}
\frac{1}{N^{1/2}} \| {\bf f}\| \leq \|{\bf f}\|_M.
\end{equation}

On the other hand, 
\begin{eqnarray*}
\|{\bf f}\|_M^4 & = & \frac{1}{N}  \sum_{k,l=1}^N
\left| \int_{a}^b \sum_{m=1}^N f_{km}(t) f_{lm}^*(t)dt \right|^2\\
         & \leq &  \sum_{k,l=1}^N  \sum_{m=1}^N \left| \int_a^b f_{km}(t)f_{lm}^*(t)dt \right|^2\\
         & \leq &   \sum_{k,l=1}^N  \sum_{m=1}^N \int_a^b |f_{km}(t)|^2dt \int_a^b |f_{lm}(t)|^2dt\\
         &\leq &  \frac{1}2  \sum_{k,l=1}^N  \sum_{m=1}^N \left( \left(\int_a^b |f_{km}(t)|^2dt\right)^2
                    + \left(\int_a^b |f_{lm}(t)|^2dt\right)^2\right)\\
         &\leq & N \sum_{k=1}^N  \sum_{m=1}^N  \left(\int_a^b |f_{km}(t)|^2dt\right)^2\\
         & \leq & N \left(  \int_a^b \sum_{k,m=1}^N |f_{km}(t)|^2dt\right)^2\\
          & = & N \|{\bf f}\|^4.
\end{eqnarray*}
This shows that
\begin{equation}\label{11}
\|{\bf f}\|_M \leq N^{1/4} \|{\bf f}\|.
\end{equation}
Combining (\ref{10}) and (\ref{11}), the equivalence (\ref{5}) with $C_1=N^{-1/2}$ and
$C_2=N^{1/4}$  between the norms $\| \cdot \|_M$ 
in ({\ref{3}) and $\|\cdot \|$ in (\ref{4}) is proved. 
{\bf q.e.d.}

Due to the equivalence of the norm $\| \cdot \|_M$ proposed in this paper 
 and the norm  $\| \cdot \|$ 
used in \cite{xia2}, all the results on orthonormal matrix-valued 
wavelets obtained in \cite{xia2} hold,  when the norm $\| \cdot \|_M$ 
for matrix-valued signals in this paper is used. 

As a remark, 
the conventional inner product for two matrices $A$ and $B$ is the scalar-valued
inner product $tr(AB^{\dagger})$ where $tr$ stands for the matrix trace.
 It is not hard to see that with this scalar-valued inner product,
the orthogonality between two matrix-valued signals, which is called 
{\em Orthogonality C},  is  the orthogonality of   two long vectors
of concatenated row vectors of two matrix-valued signals.
As mentioned in Introduction,  and  it is also not hard to see from the above 
definition, 
the orthogonality  (\ref{7}) induced from the matrix-valued inner 
product in this paper for two matrix-valued signals is the orthogonality 
between any  row vectors including 
the row vectors inside a matrix  of the two matrix-valued signals, which
is named {\em Orthogonality B} in \cite{xia2}.
Clearly Orthogonality B  is stronger 
than Orthogonality C, while it is weaker than the component-wise orthogonality
 called {\em Orthogonality A} in \cite{xia2}, commonly used in multiwavelets
\cite{strang, plonka}.

With Orthogonality A, a necessary and sufficient condition to construct 
orthonormal multiwavelets was given in \cite{plonka} that is not easy to check.
However, with Orthogonality B, an easy sufficient condition 
to construct orthonormal 
multiwavelets  was obtained in \cite{xia2}.
Furthermore, it was shown in \cite{xia2} 
that  the matrix-valued inner product (\ref{6}) and its induced Orthogonality B
provide a complete decorrelation of matrix-valued signals along time and across 
matrix components in the sense that
a complete  Karhunen-Lo\`{e}ve  expansion for matrix-valued signals
can be obtained. This may not be possible for Orthogonality A or Orthogonality C 
induced from a scalar-valued inner product \cite{Kelly, Trees}. 
In other words, the matrix-valued inner product (\ref{6}) is fundamental to 
study matrix-valued signals.

\section{Linear Independence and Gram-Schmidt Orthonormalization}\label{sec3}

Let us first introduce 
degenerate and linearly independent matrix-valued signals,
 and study their properties.

\subsection{Degenerate and Linearly Independent Matrix-Valued  Signals}

A matrix-valued signal ${\bf f}$ in $L^2(a,b; {\mathbb{C}}^{N\times N})$ 
is called {\em degenerate signal} if 
$\langle {\bf f}, {\bf f}\rangle$ does not have full rank, 
otherwise it is called {\em nondegenerate signal}. 
A sequence of  matrix-valued signals ${\bf f}_k$ in $L^2(a,b; {\mathbb{C}}^{N\times N})$, $k=1, 2, ..., K$, are  called {\em linearly independent} if the following condition holds:  if 
\begin{equation}\label{3.1}
\sum_{k=1}^K F_k {\bf f}_k \stackrel{\Delta}{=} {\bf f}
\end{equation}
for constant matrices $F_k\in \mathbb{C}^{N\times N}$, $k=1,2,...,K$, 
is degenerate,  
then the null space of  matrix $F_k^{\dagger}$ includes 
the  null space of matrix $\langle {\bf f}, {\bf f}\rangle$ 
for every $k$, $k=1,2,...,K$. 
Clearly, the above linear independence returns to the conventional one when 
all the above matrices including both $F_k$ and ${\bf f}_k$ are diagonal. 
Furthermore, if ${\bf f}=0$ in (\ref{3.1}), the above condition implies 
that all $F_k=0$, $k=1,2,...,K$, since in this case, the null space 
of $\langle {\bf f}, {\bf f}\rangle$ is the whole space 
$\mathbb{C}^{N\times N}$. This  concides with the condition of 
the conventional linear independence. 



\begin{PP}
\label{p1}
If matrix-valued signals ${\bf f}_k$, $k=1,2,..., K$, are linearly independent,
then, all signals ${\bf f}_k$, $k=1,2,...,K$,  are  nondegenerate.
\end{PP}

{\bf Proof}: Without loss of generality, assume ${\bf f}_1$ is degenerate.
Let $F_1=I_N$ and $F_k=0$ for $k=2,3,...,K$. 
Then, we have  that 
$$
\sum_{k=1}^K F_k {\bf f}_k = {\bf f}_1
$$
is degenerate, while the null space of $F_1^{\dagger}$  is $0$ only and does not include the null space  of $\langle {\bf f}_1, {\bf f}_1\rangle$. 
In other words, ${\bf f}_k$, $k=1,2,...,K$, are not linearly independent.
This contradicts the assumption in the 
proposition and therefore the proposition is proved. 
{\bf q.e.d.}

As one can see, the above concept of degenerate signal 
is similar to that of $0$ 
in the conventional linear dependence or independence.

\begin{PP}
\label{p110}
Let $G_k\in \mathbb{C}^{N\times N}$, $k=1,2,..., K$,
be $K$ constant matrices  and at least one of them have full rank. 
If matrix-valued signals ${\bf f}_k$, $k=1,2,..., K$, 
are linearly independent, then $\sum_{k=1}^K G_k {\bf f}_k$ is nondegenerate.
\end{PP}

{\bf Proof}: Without loss of generality, let us assume 
$G_1$ has full rank. 
 If $\sum_{k=1}^K G_k {\bf f}_k={\bf g}$ is degenerate,
then by the linear independence of 
${\bf f}_k$, $k=1,2,..., K$,
the null space of $G_1^{\dagger}$ cannot only contain $0$, which 
contradicts the assumption that $G_1$ has full rank. 
{\bf q.e.d.}

It is clear to see that Proposition \ref{p1} is a special case of Proposition 
\ref{p110}.

\begin{PP}
\label{p11}
If matrix-valued signals ${\bf f}_k$, $k=1,2,..., K$, are linearly independent,
then for any full rank constant matrices 
$G_k\in \mathbb{C}^{N\times N}$, $k=1,2,...,K$, matrix-valued 
signals ${\bf g}_k\stackrel{\Delta}{=}G_k {\bf f}_k$, $k=1,2,...,K$,
 are also linearly independent.
\end{PP}

{\bf Proof}: For any constant matrices $F_k \in \mathbb{C}^{N\times N}$, $k=1,2,...,K$, if 
$$
\sum_{k=1}^K F_k {\bf g}_k= \sum_{k=1}^K F_k G_k{\bf f}_k={\bf f}
$$
is degenerate, then for each $k$, $1\leq k\leq K$,
the null space of matrix $(F_kG_k)^{\dagger}=G_k^{\dagger}F_k^{\dagger}$
 includes the null space of 
matrix $\langle {\bf f}, {\bf f}\rangle $, 
since ${\bf f}_k$, $k=1,2,...,K$, are linearly independent.
Because all matrices $G_k$, $k=1,2,...,K$, have full rank, for each $k$,
$1\leq k\leq K$, the null spaces of $F_k^{\dagger}$ 
and $G_k^{\dagger}F_k^{\dagger}$ are the same, thus, 
the null space of  $F_k^{\dagger}$ 
includes the null space of $\langle {\bf f}, {\bf f}\rangle$  as well. This proves the proposition.
{\bf q.e.d.}

Similar to the conventional linear dependence of vectors, 
we have the following result for matrix-valued signals.

\begin{PP}
\label{p34}
For a matrix-valued signal ${\bf f}\in L^2(a,b; {\mathbb{C}}^{N\times N})$
and two constant matrices $A, B\in \mathbb{C}^{N\times N}$,
matrix-valued signals $A{\bf f}$ and $B{\bf f}$ are linearly dependent.
\end{PP}

{\bf Proof}: If $A{\bf f}$ and $B{\bf f}$ are linearly independent, 
then, from Proposition \ref{p1} it is easy to see 
that matrices $A$ and $B$  all have  full rank and  ${\bf f}$ is nondegenerate.
Then, we have 
$$
BA^{-1}A {\bf f}-B {\bf f}=0,
$$
which contradicts with the assumption of the linear independence
of  $A{\bf f}$ and $B{\bf f}$.
This proves the proposition.
{\bf q.e.d.}

Although it is obvious for the conventional vectors,
the result in Proposition \ref{p34} for matrix-valued signals 
may not be so, 
due to the matrix-valued coefficient multiplications as it can be seen 
from the above proof. 
We next consider more general 
linear combinations of linearly independent matrix-valued 
signals. 

For $1\leq p\leq K$, let $S_1,...,S_p$ be a partition of the index set 
$\{1,2,...,K\}$ and each $S_i$ has $K_i$ elements, where
$S_{i_1}\cap S_{i_2}=\emptyset$ for $1\leq  i_1\neq i_2\leq p$,
$\cup_{i=1}^p S_i=\{1,2,...,K\}$, 
and $1\leq K_1, ...,K_p\leq K$ with $K_1+K_2+\cdots +K_p=K$. 

\begin{PP}
\label{p222}
For each $i$, $1\leq i\leq p$,
 let $G_{k_i}\in \mathbb{C}^{N\times N}$, $k_i\in S_i$, be 
$K_i$ constant matrices and at least one of them have full rank.
If matrix-valued signals ${\bf f}_k$, $k=1,2,..., K$, are linearly independent,
then the following $p$ matrix-valued signals:
$$
\sum_{k_i\in S_i} G_{k_i} {\bf f}_{k_i},\,\,\mbox{ for  }
i=1,2,...,p,
$$
are linearly independent. 
\end{PP}

{\bf Proof}: Let $F_i\in \mathbb{C}^{N\times N}$, $i=1,2,...,p$,
 be constant matrices. Assume that 
$$
\sum_{i=1}^p F_i  \sum_{k_i\in S_i} G_{k_i} {\bf f}_{k_i} ={\bf g}
$$
is degenerate. 
Then,
$$
\sum_{i=1}^p \sum_{k_i\in S_i} F_i G_{k_i} {\bf f}_{k_i} ={\bf g},
$$
and by the linear independence of ${\bf f}_k$, $k=1,2,..., K$,
we know that the null space of $(F_iG_{k_i})^{\dagger}$ for every 
$k_i\in S_i$ and every  $i=1,...,p$ 
 contains the null space of matrix $\langle {\bf g}, {\bf g}\rangle$.
From the condition in the proposition, without loss of generality,
we may  assume that $G_{k_{i,1}}$, for some $k_{i,1}\in S_i$,
 has full rank for $1\leq i\leq p$. 
Thus,  the null space of $(F_i G_{k_{i,1}})^{\dagger}$, 
or $G_{k_{i,1}}^{\dagger}F_i^{\dagger}$, 
contains the null space of $\langle {\bf g}, {\bf g}\rangle$
for $1\leq i\leq p$. 
Since  $G_{k_{i,1}}$ 
 has full rank, the null space of $F_i^{\dagger}$
must contain the null space of $\langle {\bf g}, {\bf g}\rangle$
for $1\leq i\leq p$.
This proves the proposition.
{\bf q.e.d.}

Note that when $p=1$ in Proposition \ref{p222}, it returns 
to Proposition \ref{p110}, and when $p=K$ in Proposition \ref{p222}, it returns
to Proposition \ref{p11}. 



\begin{PP}
\label{p2}
If  ${\bf f}_k$, $k=1,2,...,K$,  form an orthonormal set 
 in $L^2(a,b; {\mathbb{C}}^{N\times N})$, then, they 
must be linearly independent.
\end{PP}

{\bf Proof}: For constant matrices $F_k\in {\mathbb C}^{N\times N}$, $k=1,2,...,K$, let
$$
\sum_{k=1}^K F_k {\bf f}_k = {\bf f}.
$$
Then, from the Parseval equality (\ref{9}), we have
$$
\sum_{k=1}^K F_k F_k^{\dagger}=\langle {\bf f}, {\bf f}\rangle.
$$
Assume that for some vector 
$u\neq 0$, we have $\langle {\bf f}, {\bf f}\rangle u=0$
but 
$F_{k_0}^{\dagger}u\neq 0$ for some $k_0$, $1\leq k_0\leq K$.
Then, 
$$
0< u^{\dagger} F_{k_0} F_{k_0}^{\dagger} u \leq 
 \sum_{k=1}^K  u^{\dagger} F_k F_k^{\dagger} u = u^{\dagger}
 \langle {\bf f}, {\bf f}\rangle u =0,
$$
which leads to a contradiction. Thus, for each $k$, $1\leq k\leq K$,
the null space of $F_k^{\dagger}$ includes
the null space of $ \langle {\bf f}, {\bf f}\rangle$.
This proves 
the linear independence of ${\bf f}_k$, $k=1,2,...,K$.
{\bf q.e.d.}

\begin{CC}
\label{c1}
Assume  ${\bf f}_k$, $k=1,2,...,K$, are  nondegenerate matrix-valued signals
and form an orthogonal set 
in $L^2(a,b; {\mathbb{C}}^{N\times N})$.
Then, ${\bf g}_k\stackrel{\Delta}{=} \langle {\bf f}_k, {\bf f}_k \rangle^{-1/2}{\bf f}_k$, $k=1,2,...,K$,
form an orthonormal set in  $L^2(a,b; {\mathbb{C}}^{N\times N})$, and 
 ${\bf f}_k$, $k=1,2,...,K$,  are linearly independent.
\end{CC}

{\bf Proof}: Since all matrix-valued signals ${\bf f}_k$ 
are nondegenerate, matrices $\langle {\bf f}_k, {\bf f}_k \rangle$
all have full rank. From the property (iv) for the matrix-valued inner
product and $\{{\bf f}_k\}$ is an orthogonal set,  for every $k,l$, $1\leq k,l \leq K$, 
$$
\langle {\bf g}_k, {\bf g}_l \rangle = 
\langle {\bf f}_k, {\bf f}_k \rangle^{-1/2}  
\langle {\bf f}_k, {\bf f}_l \rangle  \langle 
{\bf f}_l, {\bf f}_l \rangle^{-1/2} = \delta(k-l)I_N.
$$
Thus, ${\bf g}_k$, $k=1,2,...,K$, form an orthonormal
set in $L^2(a,b; {\mathbb{C}}^{N\times N})$.

Then, the linear independence of 
${\bf f}_k = \langle {\bf f}_k, {\bf f}_k\rangle^{1/2}{\bf g}_k$, 
$k=1,2,...,K$, immediately follows from 
 Propositions \ref{p11}  and \ref{p2}.
{\bf q.e.d.}

The result in Corollary \ref{c1} is consistent with 
the conventional one for vectors, i.e., any orthogonal set of nonzero 
vectors must be linearly independent. 
However, there is a difference.
In the above relationship between
orthogonality and linear independence,
matrix-valued signals need to be nondegenerate. 
Note that it is possible that a matrix-valued signal ${\bf f}$
 in an orthogonal set
 in $L^2(a,b; {\mathbb{C}}^{N\times N})$ is degenerate, 
 i.e., $\langle {\bf f}, {\bf f}\rangle$ may not necessarily have 
full rank, even though ${\bf f}\neq 0$. 
Thus, a general orthogonal set of matrix-valued signals
may not have to  be linearly independent.
This does not occur for any  orthogonal
set of nonzero signals  when a scalar-valued inner product is used.

\subsection{Gram-Schmidt Orthonormalization}

We are now ready to present the Gram-Schmidt orthonormalization for
a finite sequence of linearly independent matrix-valued signals.
Let ${\bf f}_k\in L^2(a,b; {\mathbb{C}}^{N\times N})$, $k=1,2,...,K$,
be linearly independent. The Gram-Schmidt orthonormalization for this 
sequence is as follows, which is similar to, 
but not a straightforward extension of, the conventional one, 
 due to the noncommutativity of matrix multiplications.

Since ${\bf f}_k\in L^2(a,b; {\mathbb{C}}^{N\times N})$, $k=1,2,...,K$,
are linearly independent, by Proposition  \ref{p1}, 
${\bf f}_1$ is nondegenerate,
i.e., matrix $\langle {\bf f}_1, {\bf f}_1 \rangle$ is invertible and positive
definite. Let
\begin{equation}\label{3.2}
{\bf g}_1 =\langle {\bf f}_1, {\bf f}_1\rangle^{-1/2} {\bf f}_1.
\end{equation}
Then, we have
\begin{eqnarray}
\langle {\bf g}_1, {\bf g}_1 \rangle 
 & = & \langle {\bf f}_1, {\bf f}_1\rangle^{-1/2} \int_{a}^{b}
{\bf f}_1(t) {\bf f}_1^{\dagger}(t) dt 
\langle {\bf f}_1, {\bf f}_1\rangle^{-1/2}  \nonumber \\
 & = & \langle {\bf f}_1, {\bf f}_1\rangle^{-1/2}
\langle {\bf f}_1, {\bf f}_1 \rangle 
\langle {\bf f}_1, {\bf f}_1\rangle^{-1/2}
=I_N. \label{3.3}
\end{eqnarray}

Let 
\begin{eqnarray}
 \hat{\bf g}_2 & = & {\bf f}_2 -\langle {\bf f_2}, {\bf g}_1 \rangle 
 {\bf g}_1, \label{3.4}\\
{\bf g}_2 & = &  \langle \hat{\bf g}_2, \hat{\bf g}_2\rangle^{-1/2} \hat{\bf g}_2.
\label{3.5}
\end{eqnarray}
For (\ref{3.5}) to be vaild, we need to show that 
$\hat{\bf g}_2$ in (\ref{3.4}) is nondegenerate.
In fact, if $\hat{\bf g}_2$  is degenerate, then, from 
(\ref{3.4}) and (\ref{3.2}), we have  
\begin{eqnarray*}
\hat{\bf g}_2 & = &  {\bf f}_2 - \langle {\bf f}_2, {\bf g}_1\rangle
\langle {\bf f}_1, {\bf f}_1 \rangle^{-1/2} {\bf f}_1\\
   & = & F_1 {\bf f}_1 +F_2 {\bf f}_2,
\end{eqnarray*}
where $F_2=I_N$ and $F_1=- \langle {\bf f}_2, {\bf g}_1\rangle
\langle {\bf f}_1, {\bf f}_1 \rangle^{-1/2}$. Similar to the proof 
of Proposition \ref{p1}, this contradicts
 the assumption that ${\bf f}_1$ and ${\bf f}_2$ are linearly independent. 
Therefore, it  proves that $\hat{\bf g}_2$ in (\ref{3.4}) is nondegenerate
and (\ref{3.5}) is well-defined.

Let us then check the orthogonality between ${\bf g}_1$ and $\hat{\bf g}_2$. 
From (\ref{3.4}) and (\ref{3.3}), we have
\begin{eqnarray*}
\langle \hat{\bf g}_2, {\bf g}_1\rangle
 & = &  \langle {\bf f}_2, {\bf g}_1 \rangle
-\langle {\bf f}_2, {\bf g}_1\rangle 
\langle {\bf g}_1, {\bf g}_1 \rangle\\
 & = & \langle {\bf f}_2, {\bf g}_1 \rangle
-\langle {\bf f}_2, {\bf g}_1\rangle =0.
\end{eqnarray*}
From (\ref{3.5}) and (\ref{3.3}), we have 
that ${\bf g}_1$ and ${\bf g}_2$ form an orthonormal set.

Repeat the above process and for a general $k$, $2\leq k\leq K$, we let
\begin{eqnarray}
 \hat{\bf g}_k & = & {\bf f}_k -
\sum_{l=1}^{k-1}  \langle {\bf f}_k, {\bf g}_l \rangle 
 {\bf g}_l, \label{3.6}\\
{\bf g}_k & = &  \langle \hat{\bf g}_k,  \hat{\bf g}_k\rangle^{-1/2} \hat{\bf g}_k.
\label{3.7}
\end{eqnarray}
With the same proof as the above ${\bf g}_1$ and ${\bf g}_2$, we have
the following proposition.

\begin{PP}
\label{p3}
For a linearly independent set of matrix-valued 
signals ${\bf f}_k$, $k=1,2,...,K$, 
let ${\bf g}_1, {\bf g}_2, ..., {\bf g}_K$ 
be constructed in (\ref{3.2}) and (\ref{3.4})-(\ref{3.7}).
Then, ${\bf g}_1, {\bf g}_2, ..., {\bf g}_K$ form an orthonormal set.
\end{PP}

As we can see, although the above Gram-Schmidt orthonormalization
 procedure for matrix-valued signals
is similar to the conventional one for vectors, it is not a straightforward
generalization due to 1) the noncommuntativity of matrix mulitpications
and 2) the matrix-valued inner product used in the above procedure.

We also want to make a comment on the nondegenerate and 
linear independence for matrix-valued signals.
The condition for nondegenerate matrix-valued signals is 
a weak condition. Unless the row vectors  of functions are
linearly dependent in the conventional sense, otherwise, a matrix-valued 
signal is usually nondegenerate. This can be seen from the following
result.

\begin{PP}\label{P9}
  A matrix-valued signal ${\bf f}\in L^2(a,b; {\mathbb{C}}^{N\times N})$
  is degenerate if and only if its row vectors of functions in
  $L^2(a,b)$ are linearly dependent.
\end{PP}

{\bf Proof}:  A matrix-valued signal ${\bf f}$ is degenerate if and only if 
$\langle {\bf f}, {\bf f}\rangle$ does not have full rank. Since
constant matrix $\langle {\bf f}, {\bf f}\rangle$ is Hermitian, it is true if and only if
there is a vector $u\neq 0 \in \mathbb{C}^N$ such that
$\langle u, \langle {\bf f}, {\bf f}\rangle u\rangle =0$, i.e.,
$$
\int_a^b \langle u, {\bf f}(t) {\bf f}^{\dagger}(t)u \rangle dt
=\int_a^b \| {\bf f}^{\dagger}(t)u \|_F^2 dt =0,
$$
which holds if and only if ${\bf f}^{\dagger}(t)u=0$ or
 $u^{\dagger}{\bf f}(t)=0$ almost everywhere
  for $t\in (a,b)$, i.e., the row vectors of functions of  ${\bf f}(t)$
  are linear dependent.
{\bf q.e.d.}

Thus, for a finite set of nondegenerate matrix-valued signals, they usually satisfy 
the condition for linear independence for matrix-valued signals 
defined above, i.e., they are usually 
linearly independent and therefore, they can be made to an orthonormal set
by using the above Gram-Schmidt procedure. 

Another comment on the linear independence for matrix-valued signals
is that the definition in (\ref{3.1}) is only for left multiplication 
of constant matrices $F_k$ to matrix-valued signals ${\bf f}_k$. Similar definitions for 
linear independences of matrix-valued signals with right 
constant matrix multiplications and/or mixed left and 
right constant matrix  multiplications may be possible. 
Although what is studied in this paper is for continuous-time 
matrix-valued signals, 
it can be easily generalized to discrete-time matrix-valued 
 signals (sequences of finite or infinite length).

\section{Matrix-Valued Lattices}\label{sec7}
In this section, based on the matrix-valued signal space
 with the matrix-valued inner product, we introduce matrix-valued lattices.

We first introduce matrix-valued lattices. For convenience, in what follows
we only consider 
 the Frobenius norm for matrices, i.e., 
$\| \cdot \|_M=\| \cdot \|_F$, and 
real matrix-valued signal spaces, 
and also let ${\cal R}$ denote the real matrix-valued
signal space ${\cal R}\stackrel{\Delta}{=}L^2(a,b; {\mathbb{R}}^{N\times N})$. 
Let $\mathbb{Z}^{N\times N}$ denote all $N\times N$ matrices of integer 
entries.

For a finitely 
many linearly independent real matrix-valued signals ${\bf f}_k\in {\cal R}$, $k=1,2,...,K$, let ${\cal R}^K$ denote the matrix-valued signal space linearly expanded by them, i.e.,
\begin{equation}\label{7.1}
{\cal R}^K= \left\{ \sum_{k=1}^K F_k {\bf f}_k: \,\, F_k\in \mathbb{R}^{N\times N},\,k=1,2,...,K\right\}.
\end{equation}
From what was studied in the previous section,
clearly, ${\bf f}_k$, $k=1,2,...,K$, form a basis in ${\cal R}^K$.
 The {\em matrix-valued lattice} formed by
this basis in  ${\cal R}^K$ is defined as
\begin{equation}\label{7.2}
{\cal L}=\left\{\sum_{k=1}^K F_k {\bf f}_k : \,\, F_k \in \mathbb{Z}^{N\times N}, \,
k=1,2,...,K\right\},
\end{equation}
which is a subset/subgroup of  ${\cal R}^K$.
The basis ${\bf f}_k$, $k=1,2,...,K$, is called a basis 
for the $K$ dimensional matrix-valued lattice ${\cal L}$. 



The fundamental region of this lattice 
${\cal L}$ can be defined similar to the conventional lattice as follows.
A set ${\cal F}\subset {\cal R}^K$ is called a fundamental region, if 
its translations ${\bf x}+{\cal F}=\{ {\bf x}+ {\bf f}: \,\, {\bf f}\in {\cal F}\}$ for ${\bf x}\in {\cal L}$ 
form a partition of  ${\cal R}^K$. Since the basis elements ${\bf f}_k$, $k=1,2,..., K$, are not constant real vectors as in the conventional lattices, it
would not be convenient to define the determinant of the lattice. However,
with the Gram-Schmidt orthonormalization developed in the previous 
section, we may define the determinant of the lattice 
directly as
\begin{equation}\label{7.3}
\det({\cal L})=\prod_{k=1}^K \|\hat{\bf f}_k\|_F,
\end{equation}
where $\hat{\bf f}_k$, $k=1,2,...,K$,  are from the following 
 Gram-Schmidt orthogonalization of ${\bf f}_k$, $k=1,2,...,K$,
which is  from the Gram-Schmidt 
orthonormalization in the previous section: 
$$
\hat{\bf f}_1={\bf f}_1,
$$
\begin{equation}\label{7.4}
\hat{\bf f}_k=\hat{\bf g}_k={\bf f}_k-\sum_{l=1}^{k-1} \mu_{l,k} \hat{\bf f}_l,
\end{equation}
where 
\begin{equation}\label{7.5}
\mu_{l,k}=\langle {\bf f}_k, \hat{\bf f}_l\rangle \langle \hat{\bf f}_l, 
\hat{\bf f}_l
\rangle^{-1}, \,\, l=1,2,...,k-1 \mbox{ and }k=2,3,...,K. 
\end{equation}
It is clear to see that the spaces linearly spanned by 
$\{{\bf f}_1, {\bf f}_2, ..., {\bf f}_K\}$ and 
$\{\hat{\bf f}_1, \hat{\bf f}_2, ..., \hat{\bf f}_K\}$ are the same, 
i.e., ${\cal R}^K$ in (\ref{7.1}), 
since they can be linearly (over $\mathbb{R}^{N\times N}$) represented by each other similar to the 
conventional vectors. 

From the Gram-Schmidt orthogonalization (\ref{7.4}), we have
\begin{equation}\label{7.6}
\langle {\bf f}_k,  {\bf f}_k\rangle
=\langle \hat{\bf f}_k,  \hat{\bf f}_k\rangle
+ \sum_{l=1}^{k-1} \mu_{l,k} \langle \hat{\bf f}_l,  \hat{\bf f}_l\rangle
 \mu_{l,k}^{\dagger},
\end{equation}
for $k=1,2,...,K$. 
Using Property (iii) in Section \ref{sec2},
 the identity (\ref{7.6}) implies 
\begin{equation}\label{7.7}
\|{\bf f}_k\|_F^2\leq \|\hat{\bf f}_k\|_F^2 
+\sum_{l=1}^{k-1} \|\mu_{l,k}\|_F^2 \cdot \|\hat{\bf f}_l\|_F^2,
\end{equation} 
for $k=1,2,...,K$.
From (\ref{7.6}), it is also clear that
$\|{\bf f}_k\|_F\geq \|\hat{\bf f}_k\|_F$, $k=2,...,K$. 

As a remark, we know that the 
conventional Gram-Schmidt orthogonalization
plays an important role in the LLL algorithm for 
the conventional lattice basis reduction \cite{LLL}. 
It is, however, not clear how the Gram-Schmidt 
orthogonalization for matrix-valued signals introduced 
above can be applied in matrix-valued lattice basis reduction.

\section{Conclusion}\label{sec8}

In this paper, we re-studied the matrix-valued inner product 
previously proposed to construct orthonormal matrix-valued
wavelets for matrix-valued  signal analysis \cite{xia1, xia2}
where not much on the matrix-valued inner product or 
its induced Orthogonality B, which is more fundamental, was studied. 
In order to study more on the matrix-valued inner product
and its induced  Orthogonality B, 
we first proposed a new norm for matrix-valued signals,  which is more 
consistent with the matrix-valued inner product 
than that used in \cite{xia2}, and is 
similar to that with the conventional scalar-valued inner product. 
We showed that these two norms are equivalent, which means 
that with the newly proposed norm, all the results for 
contructing orthonormal 
matrix-valued wavelets obtained in \cite{xia2} still hold. 
We then proposed the concepts of degenerate and nondegenerate matrix-valued 
signals and defined the linear independence for matrix-valued
signals,  which is different from but similar to 
the conventional linear independence for vectors. 
We also presented some  properties on the  linear independence
and the orthogonality. 
We then presented the Gram-Schmidt orthonormalization 
procedure for a set of linearly independent matrix-valued signals.
Although this procedure is similar to the conventional one 
for vectors, due to the noncommutativity of matrix multiplications
and the matrix-valued inner product used in the procedure, 
it is not a straightforward generalization. 
We finally defined matrix-valued lattices,
 where the newly introduced Gram-Schmidt orthogonalization might be applied. 


Since it was shown in \cite{xia2} that the matrix-valued inner product and
Orthogonality B provide a complete  Karhunen-Lo\`{e}ve  expansion
for matrix-valued signals, which a scalar-valued inner product may not do,
it is believed that what was studied in this paper 
for matrix-valued inner product for matrix-valued signal space
will have fundamental applications for high dimensional 
signal analysis in data science. 

As a final note, after this paper was written, it has been found that the
matrix-valued signal space with the  matrix-valued
inner product in this paper is related to Hilbert modules, see, for example, 
\cite{mod1}-\cite{mod4}. Interestingly, it was mentioned in \cite{mod5} that 
there does not exist any general notion of ``$C^*$-linear 
independence'' due to the existence of zero-divisors. We believe 
that the linear independence for matrix-valued signals introduced in
this paper is novel.

\end{document}